# New integral representations
# of the polylogarithm function

**Djurdje Cvijović**[1]

[1] Atomic Physics Laboratory, Vinča Institute of Nuclear Sciences,

P.O. Box 522, 11001 Belgrade, Serbia.

**Abstract.** Maximon has recently given an excellent summary of the properties of the Euler dilogarithm function and the frequently used generalizations of the dilogarithm, the most important among them being the polylogarithm function $Li_s(z)$. The polylogarithm function appears in several fields of mathematics and in many physical problems. We, by making use of elementary arguments, deduce several new integral representations of the polylogarithm $Li_s(z)$ for any complex $z$ for which $|z| < 1$. Two are valid for all complex $s$, whenever $\operatorname{Re} s > 1$. The other two involve the Bernoulli polynomials and are valid in the important special case where the parameter $s$ is an positive integer. Our earlier established results on the integral representations for the Riemann zeta function $\zeta(2n+1)$, $n \in N$, follow directly as corollaries of these representations.

## 1. Introduction

Recently, Maximon [1] has given an excellent summary of the defining equations and properties of the Euler dilogarithm function and the frequently used generalizations of the dilogarithm, the most important among them being the polylogarithm function $Li_s(z)$. These include integral representations, series expansions, linear and quadratic transformations, functional relations and numerical values for special arguments.

Motivated by this paper we have begun a systematic study of new integral representations for $Li_s(z)$ since it appears that only half a dozen of them can be found in the literature. Among known representations the following

$$Li_s(z) = \frac{z}{\Gamma(s)} \int_0^\infty \frac{t^{s-1}}{e^t - z} dt \qquad (\operatorname{Re} s > 0, z \notin (1, \infty)),$$

and

$$Li_s(z) = \frac{z}{\Gamma(s)} \int_1^\infty \frac{\log^{s-1} t}{t(t-z)} dt = \frac{z}{\Gamma(s)} \int_0^1 \frac{\log^{s-1}(1/t)}{1 - zt} dt \qquad (\operatorname{Re} s > 0),$$

are best known [2, pp. 30-31].

Here, by making use of fairly elementary arguments, we deduce several new integral representations of the polylogarithm function for any complex $z$ for which $|z| < 1$. Two are valid for all complex $s$, whenever $\operatorname{Re} s > 1$ (see Theorem). The other two given in Corollary 1 are valid in the important special case where the parameter $s$ is an positive integer. Our earlier published results [3] on the integral representations for the Riemann zeta function $\zeta(2n+1)$, $n \in N$, directly follow as corollaries of these representations (see Corollary 2).

## 2. The polylogarithm function

The polylogarithm function (also known as Jonquiére's function), $Li_s(z)$, is defined for any complex $s$ and $z$ as the analytic continuation of the Dirichlet series

$$Li_s(z) = \sum_{k=1}^{\infty} \frac{z^k}{k^s} \qquad (|z| < 1; s \in C) \qquad (1)$$

which is absolutely convergent for all $s$ and all $z$ inside the unit disc in the complex $z$-plane. Sometimes, $Li_s(z)$ is referred to as the polylogarithm of order $s$ and argument $z$, most frequently, however, it is simply called the polylogarithm.

The function $Li_s(z)$, for fixed $s$, has no poles or essential singularities anywhere in the finite complex $z$-plane, and, for fixed $z$, has no poles or essential singularities anywhere in the finite complex $s$-plane. It has only one essential singular point at $s = \infty$. For fixed $s$, the function has two branch points: $z = 1$ and $z = \infty$. The principal branch of the polylogarithm is chosen to be that for which $Li_s(z)$ is real for $z$ real, $0 \leq z \leq 1$, (when $s$ is real) and is continuous except on the positive real axis, where a cut is made from $z = 1$ to $\infty$ such that the cut puts the real axis on the lower half plane of $z$. Thus, for fixed $s$ the function $Li_s(z)$ is a single-valued function in the z-plane cut along the interval $(1, \infty)$ where it is continuous from below.

In the important case where the parameter $s$ is an integer, $Li_s(z)$ will be denoted by $Li_n(z)$ (or $Li_{-n}(z)$ when negative). $Li_0(z)$ and all $Li_{-n}(z)$, $n = 1, 2, \ldots$, are rational functions in $z$. $Li_n(z)$, $n = 1, 2, \ldots$, is the polylogarithm of order $n$, (i.e. the n-th order polylogarithm). The special case $n = 1$ is the ordinary logarithm $Li_1(z) = -\log(1-z)$, while the cases $n = 2, 3, 4, \ldots,$ are classical polylogarithms known, respectively, as dilogarithm, trilogarithm, quadrilogarithm, etc.

For more details and an extensive list of references in which polylogarithms appear in physical and mathematical problems we refer the reader to Maximon [1]. The polylogarithm $Li_n(z)$ of order $n = 1, 2, 3, \ldots$ is thoroughly covered in Lewin's standard text [4], while many formulae involving $Li_s(z)$ can be found in [2, pp. 30-31] and [5]. Berndt's treatise [6] can serve as an excellent introductory text on the polylogarithm (and numerous related functions) and as an encyclopaedic source.

## 2. Statement of the results

In preparation for the statement of the results we need to make several definitions.

The Riemann zeta function and the Hurwitz zeta function, $\zeta(s)$ and $\zeta(s,a)$, are, respectively, defined by means of the series [7, p. 807, Eqs. 23.2.1, 23.2.19 and 23.2.20]:

$$\zeta(s) = \begin{cases} \sum_{k=1}^{\infty} \dfrac{1}{k^s} = \dfrac{1}{1-2^{-s}} \sum_{k=0}^{\infty} \dfrac{1}{(2k+1)^s} & (\operatorname{Re} s > 1) \\ \dfrac{1}{1-2^{1-s}} \sum_{k=1}^{\infty} \dfrac{(-1)^{k-1}}{k^s} & (\operatorname{Re} s > 0;\ s \neq 1) \end{cases} \qquad (2)$$

and [2, pp. 24-27]

$$\zeta(s,a) = \sum_{k=0}^{\infty} \frac{1}{(k+a)^s} \qquad (\operatorname{Re} s > 1,\ a \neq 0, -1, -2, \ldots). \qquad (3)$$

Both are analytic over the whole complex plane, except at $s = 1$, where they have a simple pole.

Next, for any real $x$ and any complex $s$ with $\operatorname{Re} s > 1$, we define

$$S_s(x) = \sum_{k=1}^{\infty} \frac{\sin(kx)}{k^s}, \tag{4a}$$

and

$$C_s(x) = \sum_{k=1}^{\infty} \frac{\cos(kx)}{k^s}. \tag{4b}$$

We also use the Bernoulli polynomials of degree $n$ in $x$, $B_n(x)$, defined for each nonnegative integer $n$ by means of their generating function [7, p. 804, Eq. 23.1.1]

$$\frac{te^{tx}}{e^t - 1} = \sum_{n=0}^{\infty} B_n(x) \frac{t^n}{n!} \qquad (|t| < 2\pi), \tag{5a}$$

and given explicitly in terms of the Bernoulli numbers $B_n = B_n(0)$ by [5, p. 765]

$$B_n(x) = \sum_{k=0}^{n} \binom{n}{k} B_k x^{n-k}. \tag{5b}$$

Our results are as follows.

**Theorem.** Assume that $s$ and $z$ are complex numbers, let $Li_s(z)$ be the polylogarithm function and let $S_s(x)$ and $C_s(x)$ be defined as in Eqs. (4a) and (4b), respectively. If $\operatorname{Re} s > 1$ and $|z| < 1$, then

$$Li_s(z) = \frac{1}{\delta}\int_0^\delta S_s(2\pi t)\frac{2z\sin(2\pi t)}{1-2z\cos(2\pi t)+z^2}dt, \qquad (6a)$$

$$Li_s(z) = \frac{1}{\delta}\int_0^\delta C_s(2\pi t)\frac{1-z^2}{1-2z\cos(2\pi t)+z^2}dt = \qquad (6b)$$

$$= \frac{1}{\delta}\int_0^\delta C_s(2\pi t)\frac{2z(\cos(2\pi t)-z)}{1-2z\cos(2\pi t)+z^2}dt. \qquad (6c)$$

Here, and throughout the text, we set either $\delta=1$ or $\delta=1/2$.

**Remark 1.** Observe that our integral representations in the above Theorem essentially involve either the Hurwitz zeta function or the generalized Clausen function. Indeed, between the functions $S_s(x)$ and $C_s(x)$ and the Hurwitz zeta function defined as in (3) there exists the relationship [8, p. 726, Entry 5.4.2.1]

$$\begin{Bmatrix}S_s(2\pi t)\\C_s(2\pi t)\end{Bmatrix} = \frac{(2\pi)^s}{4\Gamma(s)}\begin{Bmatrix}\csc(\pi s/2)\\\sec(\pi s/2)\end{Bmatrix}\left[\zeta(1-s,t)\mp\zeta(1-s,1-t)\right]$$

$$(0\le t\le 1; \operatorname{Re}s>1;\begin{Bmatrix}s\ne 2n\\s\ne 2n+1\end{Bmatrix}, n\in N),$$

where $\Gamma(s)$ is the familiar gamma function. In the cases when this relationship does not hold $S_s(x)$ and $C_s(x)$ are used to define the generalized Clausen function [4, p. 281, Eq. 3]

$$Cl_n(x) = \begin{Bmatrix}C_n(x) \text{ if } n \text{ is odd}\\S_n(x) \text{ if } n \text{ is even}\end{Bmatrix} \quad (n\in N).$$

What is most important, however, is that $S_{2n-1}(x)$ and $C_{2n}(x)$ ($n\in N$) are expressible in terms of the Bernoulli polynomials and this makes possible our corollaries.

**Corollary 1.** Assume that $n$ is a positive integer, let $Li_n(z)$ be the polylogarithm function and let $B_n(x)$ be the Bernoulli polynomials of degree $n$ in $x$ given by Eqs. (5). Then, provided that $|z| < 1$, we have

$$Li_{2n-1}(z) = (-1)^n \frac{1}{\delta} \frac{(2\pi)^{2n-1}}{2(2n-1)!} \int_0^\delta B_{2n-1}(t) \frac{2z\sin(2\pi t)}{1 - 2z\cos(2\pi t) + z^2} dt, \qquad (7a)$$

$$Li_{2n}(z) = (-1)^{n-1} \frac{1}{\delta} \frac{(2\pi)^{2n}}{2(2n)!} \int_0^\delta B_{2n}(t) \frac{1 - z^2}{1 - 2z\cos(2\pi t) + z^2} dt = \qquad (7b)$$

$$= (-1)^{n-1} \frac{1}{\delta} \frac{(2\pi)^{2n}}{2(2n)!} \int_0^\delta B_{2n}(t) \frac{2z(\cos(2\pi t) - z)}{1 - 2z\cos(2\pi t) + z^2} dt. \qquad (7c)$$

**Corollary 2.** Assume that $n$ is a positive integer, let $\zeta(s)$ be the Riemann zeta function defined as in (2) and let $B_n(x)$ be the Bernoulli polynomials of degree $n$ in $x$ given by Eqs. (5). We then have:

$$\zeta(2n+1) = (-1)^{n-1} \frac{1}{\delta} \frac{(2\pi)^{2n+1}}{2(2n+1)!} \int_0^\delta B_{2n+1}(t) \cot(\pi t) dt, \qquad (8a)$$

$$\zeta(2n+1) = (-1)^{n-1} \frac{2^{2n}}{2^{2n}-1} \frac{1}{\delta} \frac{(2\pi)^{2n+1}}{2(2n+1)!} \int_0^\delta B_{2n+1}(t) \tan(\pi t) dt. \qquad (8b)$$

**Remark 2.** The integral for $\zeta(2n+1)$ in (8a) when $\delta = 1$ is well known [7, p. 807, Eq. 23.2.17] while the other integrals in Corollary 2 were recently deduced [3, Theorem 1]. However, we have failed to find the integral representations given in Theorem and Corollary 1 in the literature.

**Remark 3.** The integral representations given in (6) and (7) hold for $|z| < 1$. However, our results may be extended to any complex $z$, $|z| > 1$, by means of the inversion formula [2, pp. 30-31; 5, pp. 762-763]

$$Li_s(z) = \frac{(2\pi)^s}{\Gamma(s)} e^{i\pi s/2} \zeta\left(1-s, \frac{\log(-z)}{2\pi i} + \frac{1}{2}\right) - e^{i\pi s} Li_s\left(\frac{1}{z}\right) \qquad z \notin (0,1),$$

where $\zeta(s,a)$ is the Hurwitz zeta function (3) or, by means of the following particularly simple inversion formula

$$Li_n(z) = (-1)^{n-1} Li_n\left(\frac{1}{z}\right) - \frac{(2\pi i)^n}{n!} B_n\left(\frac{\log z}{2\pi i}\right) \qquad n = 0, 1, 2, \ldots,$$

valid for the n-th order polylogarithm.

**Remark 4.** We remark that the existence of the integrals in (8) is assured since the integrands on $[0, \alpha]$, $0 < \alpha < 1$, have only removable singularities. This can be demonstrated easily by making use of some basic properties of $B_n(x)$. For instance, knowing that the odd-indexed Bernoulli numbers $B_n$, apart from $B_1 = -1/2$, are zero [7, p. 805, Eq. 23.1.19], we have

$$\lim_{t \to 1/2} B_{2n+1}(t) \tan(\pi t) = \lim_{t \to 1/2} \frac{B_{2n+1}(t)}{\cos(\pi t)} = \lim_{t \to 1/2} \frac{(2n+1)B_{2n}(t)}{-\pi \sin(\pi t)} = (1 - 2^{1-2n})(2n+1)B_{2n}/\pi$$

since $B_n(1/2) = (2^{1-n} - 1)B_n$ [7, p. 805, Eq. 23.1.21]

### 3. Proof of the results

In what follows we shall need the following lemma and we provide two proofs for it. The second proof was suggested by one of the anonymous referees and it makes use of well-known Chebyshev polynomials and Fourier series.

**Lemma.** Assume that $n$ is a positive integer and that $\delta = 1$ and $\delta = \frac{1}{2}$. Then we have:

$$\int_0^\delta \begin{Bmatrix} \cos(n2\pi t) \\ \sin(n2\pi t) \end{Bmatrix} \frac{2z\sin(2\pi t)}{1 - 2z\cos(2\pi t) + z^2} dt = \begin{Bmatrix} 0 \\ \delta z^n \end{Bmatrix}, \qquad (9a)$$

$$\int_0^\delta \begin{Bmatrix} \cos(n2\pi t) \\ \sin(n2\pi t) \end{Bmatrix} \frac{1 - z^2}{1 - 2z\cos(2\pi t) + z^2} dt = \begin{Bmatrix} \delta z^n \\ 0 \end{Bmatrix}. \qquad (9b)$$

**First Proof of Lemma.** It is clear that the integrals in (9) can be rewritten as follows:

$$\int_0^{2\delta\pi} \begin{Bmatrix} \cos(nt) \\ \sin(nt) \end{Bmatrix} \frac{2z\sin t}{1 - 2z\cos t + z^2} dt = \begin{Bmatrix} 0 \\ 2\delta\pi z^n \end{Bmatrix}, \qquad (10a)$$

$$\int_0^{2\delta\pi} \begin{Bmatrix} \cos(nt) \\ \sin(nt) \end{Bmatrix} \frac{1 - z^2}{1 - 2z\cos t + z^2} dt = \begin{Bmatrix} 2\delta\pi z^n \\ 0 \end{Bmatrix}. \qquad (10b)$$

We first establish the case $\delta = 1$ in (10). For any positive integer $n$ and arbitrary complex $z$ consider the integrals $I_1$ and $I_2$ with parameters given by

$$I_1(n, z) = \int_0^{2\pi} e^{int} \frac{2z\sin t}{1 - 2z\cos t + z^2} dt =$$

$$= \int_0^{2\pi} \cos(nt) \frac{2z\sin t}{1 - 2z\cos t + z^2} dt + i \int_0^{2\pi} \sin(nt) \frac{2z\sin t}{1 - 2z\cos t + z^2} dt, \qquad (11a)$$

and by

$$I_2(n, z) = \int_0^{2\pi} e^{int} \frac{1 - z^2}{1 - 2z\cos t + z^2} dt =$$

$$= \int_0^{2\pi} \cos(nt) \frac{1-z^2}{1-2z\cos t + z^2} dt + i \int_0^{2\pi} \sin(nt) \frac{1-z^2}{1-2z\cos t + z^2} dt. \tag{11b}$$

In order to derive the integrals $I_1$ and $I_2$ we make use of contour integration and calculus of residues. By setting

$$\tau = e^{it}, \quad \cos t = \frac{1}{2}(\tau + 1/\tau) \quad \text{and} \quad \sin t = \frac{1}{2i}(\tau - 1/\tau)$$

where $i = \sqrt{-1}$, we arrive at

$$I_1(n,z) = \oint_{|\tau|=1} \frac{i\tau^n(\tau^2-1)}{(\tau-z)(\tau-1/z)} \frac{d\tau}{i\tau} = \oint_{|\tau|=1} \frac{\tau^{n-1}(\tau^2-1)}{(\tau-z)(\tau-1/z)} d\tau =$$

$$= 2\pi i \operatorname*{Res}_{\tau=z} \frac{\tau^{n-1}(\tau^2-1)}{(\tau-z)(\tau-1/z)} = 2\pi i\, z^n, \tag{12a}$$

and

$$I_2(n,z) = \oint_{|\tau|=1} \frac{\tau^{n+1}(z^2-1)}{z(\tau-z)(\tau-1/z)} \frac{d\tau}{i\tau} = \oint_{|\tau|=1} \frac{\tau^n(z^2-1)}{iz(\tau-z)(\tau-1/z)} d\tau =$$

$$= 2\pi i \operatorname*{Res}_{\tau=z} \frac{\tau^n(z^2-1)}{iz(\tau-z)(\tau-1/z)} = 2\pi i(-iz^n) = 2\pi z^n. \tag{12b}$$

Combining (11) and (12) together and equating the real and imaginary parts on both sides gives the integrals in (10) where $\delta = 1$. In this way we evaluate the integrals in (9) where $\delta = 1$. Observe that in both integrals in (12) the contour is the unit circle and is traversed in the positive (counterclockwise) direction and the only singularities of the integrands that lie inside the contour are at $\tau = z$.

Next, it can be easily shown that the case $\delta = 1/2$ in (10) reduces to the above considered case $\delta = 1$. Indeed, in view of the following well-known property [8, p. 272 and 273, Eqs. 2.1.2.20 and 2.1.2.21]

$$\int_0^{2a} f(t)dt = \begin{cases} 0 & f(2a-t) = -f(t) \\ 2\int_0^a f(t)dt & f(2a-t) = f(t) \end{cases}$$

it suffices to show that both integrands in (10) are such that $f(2\pi - t) = f(t)$.

This completes the proof of our lemma.

**Second Proof of Lemma.** Consider the Chebyshev polynomials of the first and second kind, $T_m(x)$ and $U_m(x)$, defined by [7, p. 776, Eqs. 22.3.15 and 22.3.16]

$$T_m(\cos\theta) = \cos(m\theta) \quad \text{and} \quad U_m(\cos\theta) = \sin[(m+1)\theta]/\sin\theta,$$

and recall that their generating functions are given by [7, p. 783, Eqs. 22.9.9 and 22.9.10]

$$\frac{1-z^2}{1-2xz+z^2} = 1 + 2\sum_{m=1}^{\infty} T_m(x) z^m, \quad \frac{1}{1-2xz+z^2} = \sum_{m=0}^{\infty} U_m(x) z^m \quad (|z| < 1).$$

Now, by applying the orthogonality relation of the trigonometric functions we, for instance, have the sought formulas in (9b)

$$\int_0^{\delta} \begin{Bmatrix} \cos(n2\pi t) \\ \sin(n2\pi t) \end{Bmatrix} \frac{1-z^2}{1-2z\cos(2\pi t)+z^2} dt = \frac{1}{2\pi} \int_0^{2\delta\pi} \begin{Bmatrix} \cos(nt) \\ \sin(nt) \end{Bmatrix} \frac{1-z^2}{1-2z\cos t+z^2} dt$$

$$= \frac{1}{2\pi} \int_0^{2\delta\pi} \begin{Bmatrix} \cos(nt) \\ \sin(nt) \end{Bmatrix} \left(1 + 2\sum_{m=1}^{\infty} T_m(\cos t) z^m\right) dt = \frac{1}{2\pi} \int_0^{2\delta\pi} \begin{Bmatrix} \cos(nt) \\ \sin(nt) \end{Bmatrix} dt$$

$$+ 2\sum_{m=1}^{\infty} z^m \frac{1}{2\pi} \int_0^{2\delta\pi} \left(\begin{Bmatrix} \cos(nt) \\ \sin(nt) \end{Bmatrix} \cos(mt)\right) dt = \begin{Bmatrix} \delta z^n \\ 0 \end{Bmatrix}.$$

The integrals given in (9a) follow in a similar manner.

**Proof of Theorem.** The proof of Theorem rests on the above Lemma. From (9a) we obtain the formula

$$z^k = \frac{1}{\delta}\int_0^{\delta} \sin(k2\pi t)\frac{2z\sin(2\pi t)}{1-2z\cos(2\pi t)+z^2}dt \qquad \delta = 1, 1/2; k = 1, 2, \ldots,.$$

Dividing the formula by $k^s$ and summing over k, we have

$$\sum_{k=1}^{\infty}\frac{z^k}{k^s} = \sum_{k=1}^{\infty}\frac{1}{\delta}\int_0^{\delta}\frac{\sin(k2\pi t)}{k^s}\frac{2z\sin(2\pi t)}{1-2z\cos(2\pi t)+z^2}dt \qquad \mathrm{Re}\, s > 1,$$

so that

$$\sum_{k=1}^{\infty}\frac{z^k}{k^s} = \frac{1}{\delta}\int_0^{\delta}\sum_{k=1}^{\infty}\frac{\sin(k2\pi t)}{k^s}\frac{2z\sin(2\pi t)}{1-2z\cos(2\pi t)+z^2}dt \qquad \mathrm{Re}\, s > 1. \qquad (13)$$

Now, the required integral formula in (6a) directly follows from the last expression in light of the definitions of $Li_s(z)$ in (1) and $S_s(x)$ in (4a). It should be noted that inverting the order of summation and integration on the right-side of (13) is justified by absolute convergence of the series involved.

Starting from (9b) the formula in (6b) is derived in precisely the same way. In order to prove (6c) note that

$$\frac{1-z^2}{1-2z\cos(2\pi t)+z^2} = 1 + \frac{2z(\cos(2\pi t)-z)}{1-2z\cos(2\pi t)+z^2},$$

thus we have

$$\int_0^\delta C_s(2\pi t)\frac{1-z^2}{1-2z\cos(2\pi t)+z^2}dt = \int_0^\delta C_s(2\pi t)\frac{2z(\cos(2\pi t)-z)}{1-2z\cos(2\pi t)+z^2}dt,$$

since

$$\int_0^\delta C_s(2\pi t)dt = \sum_{k=1}^\infty \frac{1}{k^s}\int_0^\delta \cos(2k\pi t)dt = 0.$$

This proves the theorem.

**Proof of Corollary 1.** We need the Fourier expansions for the Bernoulli polynomials $B_n(x)$ [7, p. 805, Eqs. 23.1.17]

$$B_{2n-1}(x) = (-1)^n \frac{2(2n-1)!}{(2\pi)^{2n-1}} \sum_{k=1}^\infty \frac{\sin(2k\pi x)}{k^{2n-1}}$$

where $0 \le x \le 1$ for $n = 2, 3, \ldots$, $0 < x < 1$ for $n = 1$, and [7, p. 805, Eqs. 23.1.18]

$$B_{2n}(x) = (-1)^{n-1} \frac{2(2n)!}{(2\pi)^{2n}} \sum_{k=1}^\infty \frac{\cos(2k\pi x)}{k^{2n}}$$

where $0 \le x \le 1$ for $n = 1, 2, \ldots$. These Fourier expansions can be rewritten as follows

$$S_{2n-1}(x) = (-1)^n \frac{(2\pi)^{2n-1}}{2(2n-1)!} B_{2n-1}(x),$$

$$C_{2n}(x) = (-1)^{n-1} \frac{(2\pi)^{2n}}{2(2n)!} B_{2n}(x),$$

in terms of the functions $S_s(x)$ and $C_s(x)$ defined in (4a) and (4b).

Finally, the integral formulae for n-th polylogarithms proposed in (7) are obtained from the above expansions in conjunction with Theorem.

**Proof of Corollary 2.** First, note that for $z = 1$ and $z = -1$ the polylogarithm reduces to the Riemann zeta function

$$Li_s(1) = \zeta(s) \qquad (\operatorname{Re} s > 1), \qquad (14a)$$

$$Li_s(-1) = (2^{1-s} - 1)\zeta(s) \qquad (s \neq 1). \qquad (14b)$$

Secondly, it is not difficult to show that

$$\lim_{z \to 1} \frac{2z \sin(2\pi t)}{1 - 2z \cos(2\pi t) + z^2} = \cot(\pi t) \qquad (15a)$$

and

$$\lim_{z \to -1} \frac{2z \sin(2\pi t)}{1 - 2z \cos(2\pi t) + z^2} = \tan(\pi t). \qquad (15b)$$

Finally, we deduce the integral formulae in (8) starting from (7a) and by employing the expressions in (14) and in (15).